\documentstyle[amstex,amscd]{article}
        \catcode`\"=12
\mathchardef\za="710B  
\mathchardef\zb="710C  
\mathchardef\zg="710D  
\mathchardef\zd="710E  
\mathchardef\zve="710F 
\mathchardef\zz="7110  
\mathchardef\zh="7111  
\mathchardef\zvy="7112 
\mathchardef\zi="7113  
\mathchardef\zk="7114  
\mathchardef\zl="7115  
\mathchardef\zm="7116  
\mathchardef\zn="7117  
\mathchardef\zx="7118  
\mathchardef\zp="7119  
\mathchardef\zr="711A  
\mathchardef\zs="711B  
\mathchardef\zt="711C  
\mathchardef\zu="711D  
\mathchardef\zvf="711E 
\mathchardef\zq="711F  
\mathchardef\zc="7120  
\mathchardef\zw="7121  
\mathchardef\ze="7122  
\mathchardef\zy="7123  
\mathchardef\zvp="7124  
\mathchardef\zvr="7125 
\mathchardef\zvs="7126 
\mathchardef\zf="7127  
\mathchardef\zG="7000  
\mathchardef\zD="7001  
\mathchardef\zY="7002  
\mathchardef\zL="7003  
\mathchardef\zX="7004  
\mathchardef\zP="7005  
\mathchardef\zS="7006  
\mathchardef\zU="7007  
\mathchardef\zF="7008  
\mathchardef\zC="7009  
\mathchardef\zW="700A  

        \catcode`\"=\active
\def\*{{\textstyle *}}
\def\s*{{\scriptstyle *}}
\newfont{\ssanss}{cmss8 scaled 900}
\newfont{\smss}{cmss8 scaled 800}
\def\sss#1{\hbox{\ssanss #1}}
\newfont{\teneurm}{eurm10  scaled \magstep 1}
\def\proof{\begin{pf}}
\def\endproof{\hfill \vrule height4pt width6pt depth2pt \end{pf}}
\def\R{{\Bbb R}}
\def\Z{{\Bbb Z}}
\def\sT{\mbox{\sf T}}
\def\ssT{\sss T}

\def\ezt{\mbox{\teneurm \char'034}}
\def\ezp{\mbox{\teneurm \char'031}}
\def\ezw{\mbox{\teneurm \char'041}}
\def\xi{\operatorname{i}}
\def\xv{\operatorname{v}}
\def\xd{\operatorname{d}\!}
\def\sdt{\xd_{\hbox{\smss T}}}
\def\dt{\xd_{\sss T}}
\def\vt{\xv_{\sss T}}
\def\ll{{\pounds}}
\def\G{{\cal G}}
\def\I{{\cal J}}
\def\ad{\operatorname{ad}}
\newsymbol\leqs 1336
\newsymbol\geqs 133E
\def\Ref{\ref}
\newtheorem{defn}{Definition}[section]
\newtheorem{thm}{Theorem}[section]

\newtheorem{lem}{Lemma}[section]
\newenvironment{pf}{{\sc Proof.}\ }{}

\normalbaselines
\catcode `\@=11
\@addtoreset{equation}{section}

\def\section{\@startsection {section}{1}{\z@}{-3.5ex plus -1ex minus
     -.2ex}{2.3ex plus .2ex}{\normalsize\bf}}
\def\subsection{\@startsection{subsection}{2}{\z@}{-3.25ex plus -1ex minus
 -.2ex}{1.5ex plus .2ex}{\normalsize\bf}}

\def\thebibliography#1{\section*{References\markboth
  {REFERENCES}{REFERENCES}}\list
  {[\arabic{enumi}]}{\settowidth\labelwidth{[#1]}\leftmargin\labelwidth
  \advance\leftmargin\labelsep
  \usecounter{enumi}}
  \def\newblock{\hskip .11em plus .33em minus -.07em}
  \sloppy
  \sfcode`\.=1000\relax}
 
\catcode `\@=12

\begin{document}

        \vspace*{2.5cm} \noindent \hspace{1.2cm}{ \bf LIE ALGEBROIDS AND
POISSON-NIJENHUIS STRUCTURES} \vspace{1.3cm}
\footnote[1]{Supported by KBN, grant No 2 PO3A 074 10}\\
\noindent \hspace*{1in}
        \begin{minipage}{13cm}
        Janusz Grabowski$^{1}$ and Pawe\l\ Urba\'nski$^{2}$
\vspace{0.3cm}\\
 $^{1}$ Institute of Mathematics, Polish Academy of  Sciences,\\
      \makebox[3mm]{ } \'Sniadeckich 8, 00-950 Warszawa, Poland \\
        \makebox[3mm]{ } Institute of Mathematics, University of Warsaw, \\
        \makebox[3mm]{ } Banacha 2, 02-097 Warszawa, Poland\\
        \makebox[3mm]{ } {\sf e-mail:\ jagrab\verb+@+mimuw.edu.pl}\\
 $^{2}$ Department of Mathematical Methods in Physics,  University of \\
        \makebox[3mm]{ } Warsaw, Ho\.za 74, 00-682 Warszawa, Poland. \\
\makebox[3mm]{ } {\sf e-mail:\ urbanski\verb+@+fuw.edu.pl}
\end{minipage}

\vspace*{0.5cm}

        \begin{abstract}
        \noindent   Poisson-Nijenhuis structures  for an arbitrary
Lie algebroid are defined and studied by means of complete lifts
of tensor fields.
        \end{abstract}

\section{\hspace{-4mm}.\hspace{2mm}INTRODUCTION }

In our previous paper \cite{G-U2}, certain definitions and
constructions of graded Lie brackets and lifts of tensor fields
over a manifold were generalized to arbitrary Lie algebroids.
Since Poisson-Nijenhuis structures seem to fit very well to the
Lie algebroid language and, as it was recently  shown by
Kosmann-Schwarzbach in \cite{KS}, they give examples of Lie
bialgebroid structures in the sense of Mackenzie and Xu
\cite{M-X}, we would like to present in this note a Lie
algebroid approach to Poisson-Nijenhuis structures.

We start with the definition of a pseudo-Lie algebroid structure
on a vector bundle $E$, as a slight generalization of the notion
of a Lie algebroid and we show that such structures are
determined by special tensor fields $\zL$ on the dual bundle $E^\*$.

        Then, we define the complete lift $\dt^\zL$, which reduces
to the classical tangent  lift $\dt$ in the case of the tangent
bundle $E=\sT M$. We prove that, when we start from $P\in
\zG(M,\wedge ^2E)$, the complete lift $\dt^\zL(P)$ corresponds to
a bracket on sections of $E^\*$, which, in the classical case,
is  the Fuchssteiner-Koszul bracket on 1-forms. In the case of a
Lie algebroid over a single point, $\dt^\zL(P)$ is closely
related to the modified Yang-Baxter equation. Deforming the Lie
algebroid bracket by a (1,1) tensor $N$, we find the
corresponding bivector field $\zL_N$ on $E^\*$. Assuming some
compatibility conditions for $N$ and $P$, we can define a
Poisson-Nijenhuis structure for a Lie algebroid which  provides
a whole list of Lie bialgebroid structures.

This unified approach to Poisson and Nijenhuis structures,
including the classical case as well, as the case of a real Lie
algebra, makes possible to understand common aspects of the
theory, which were previously seen separately for different
models.

\section{\hspace{-4mm}.\hspace{2mm}TANGENT LIFTS FOR PRE-LIE ALGEBROIDS }

        Let $M$  be a manifold and let $\zt \colon E\rightarrow M$
be a vector bundle. By $\zF(\zt)$ we denote the graded exterior
algebra generated by sections of $\zt$: $\zF(\zt) = \oplus_{k\in
\Z} \zF^k(\zt)$, where $\zF^k(\zt) =\zG(M, \wedge ^kE)$ for
$k\geqs 0$ and $\zF^k(\zt) = \{ 0\}$ for $k<0$. Elements of
$\zF^0(\zt)$ are  functions on $M$, i.e., sections of the bundle
$\wedge ^0 E = M\times \R$. Similarly, by $\otimes(\zt) $ we
denote the tensor algebra $\otimes(\zt) = \oplus _{k\in \Z}
\otimes^k(\zt)$, where $\otimes^k(\zt) = \zG(M, {\otimes_M
^k}E)$. The dual vector bundle we denote by $\zp\colon E^\*
\rightarrow M$. For the tangent bundle $\ezt_M\colon \sT M
\rightarrow M$, $\zF(\ezt_M)$ is the exterior algebra of
multivector fields, and for the cotangent bundle $\ezp_M \colon
\sT^\* M\rightarrow M$, we get $\zF(\ezp_M)$, the exterior
algebra of differential forms on $M$.

        The cotangent bundle is endowed with the canonical symplectic
form $\ezw_M$ and the corresponding canonical Poisson tensor
$\zL_M$.

        \begin{defn}
        A {\em pseudo-Lie algebroid structure} on a vector bundle $\zt\colon
E\rightarrow M$  is a bracket (bilinear operation) $[\,,\,]$ on
the space $\zF^1(\zt) = \zG(M,E)$ of sections of $\zt$ and
vector bundle morphisms $\za_l,\za_r \colon E\rightarrow \sT M$
(called the left- and right-anchor, respectively), such that
        \begin{equation} \label{Fp1}
        [fX,gY] = f\za_l(X)(g)Y -g\za_r(Y)(f)X + fg[X,Y]
        \end{equation}
        for all $X,Y \in \zG(E)$ and $f,g\in C^\infty (M)$.
        \end{defn}

        A pseudo-Lie algebroid, with a skew-symmetric bracket
$[\,,\,]$ (in this case the left and right anchors coincide), is
called a {\em pre-Lie algebroid}.

        A pre-Lie algebroid is called a {\em Lie algebroid}  if the
bracket $[\,,\,]$ satisfies the Jacobi identity, i.e., if it
provides $\zF^1(\zt)$ with a Lie algebra structure.
\medskip

In the following, we establish a correspondence between
pseudo-Lie algebroid structures on $E$ and 2-contravariant
tensor fields on the  bundle manifold $E^\*$ of the dual vector
bundle $\zp \colon E^\* \rightarrow M$.
        Let $X\in \zF^1(\zt)$. We define a function $\zi_{E^\s*}X$ on
$E^\*$ by the formula
        $$ E^\* \ni a\mapsto  \zi_{E^\s*}X (a) = \langle X(\zp(a)),a \rangle ,
 $$
        where $\langle \,,\,\rangle $ is the canonical pairing
between $E$ and $E^\*$.

        Let $\zL\in \zG(E,\sT E\otimes _E \sT E)$  be a
2-contravariant tensor field on $E$.  We say that $\zL$ is {\em
linear} if, for each pair $(\zm,\zn)$ of sections of $\zp$, the
function $\zL(\xd \zi_{E}\zm, \xd \zi_{E}\zn)$,         defined on
$E$, is linear.

        For each 2-contravariant tensor $\zL$,  we define a bracket
$\{\,,\,\}_\zL$  of functions by the formula
        $$
                        \{f,g\}_\zL = \zL(\xd f, \xd g).
                                                                        $$
        \begin{thm}\label{Cp1}
                For every pseudo-Lie algebroid structure on $\zt\colon
E\rightarrow M$, with the bracket $ [\,,\,]$ and anchors
$\za_l,\za_r$, there is a unique 2-contravariant linear tensor
field $\zL$ on $E^\*$ such that
                \begin{equation}\label{Fp2}
        \zi_{E^\s*}[X,Y] = \{ \zi_{E^\s*}X, \zi_{E^\s*}Y\}_\zL
        \end{equation}
        and
        \begin{equation}\label{Fp3}
        \zp^\*\left(\za_l(X)(f)\right) = \{ \zi_{E^\s*}X, \zp^\*f\}_\zL, \ \
\zp^\*\left(\za_r(X)(f)\right)= \{ \zp^\*f, \zi_{E^\s*}X\}_\zL,
        \end{equation}
        for all $X,Y \in \zF^1(\zt)$ and $f\in C^\infty(M)$.

        Conversely, every 2-contravariant linear tensor field $\zL$
on $E^\*$ defines a pseudo-Lie algebroid on $E$ by the formulae
\ref{Fp2} and \ref{Fp3}.

        The pseudo-Lie algebroid structure on $E$ is a pre-Lie
algebroid structure (resp. a Lie algebroid structure) if and only
if the tensor $\zL$ is skew-symmetric (resp. $\zL$ is a Poisson
tensor).
        \end{thm}
        \proof
        We shall use local coordinates.
        Let $(x^a)$ be a local coordinate system on $M$ and let
$e_1,\dots, e_n$ be  a basis of local sections of $E$. We denote
by $ e^{*1},\dots, e^{*n}$ the dual basis of local sections of
$E^\*$ and by $(x^a,y^i)$ (resp. $(x^a,\zx_i)$) the corresponding
coordinate system on $E$ (resp. $E^\*$), i.e., $\zi_{E^\s*}e_i
=\zx_i$ and $\zi_{E}e^{*i}=y^i$.

        It is easy to see that every linear 2-contravariant tensor
$\zL$ on $E^\*$ is of the form
        \begin{equation}\label{Fp4}
        \zL = c^k_{ij}\zx_k \partial _{\zx_i} \otimes \partial
_{\zx_j} + \zd^a_i \partial _{\zx_i} \otimes \partial _{x^a} -
\zs^a_i \partial _{x^a}\otimes \partial _{\zx_i},
        \end{equation}
        where $c^k_{ij}, \zd^a_i \text{\ and\ } \zs^a_i$ are
functions of $x^a$.
        The correspondence between $\zL$ and a pseudo-Lie algebroid
structure is given by the formulae
        \begin{equation}\label{Fp5}
        \begin{split}
          [e_i,e_j]= [e_i,e_j]^\zL  &=c_{ij}^ke_k\\
        \za^\zL_l(e_i) &= \zd^a_i \partial _{x^a}\\
        \za^\zL_r(e_i) &= \zs^a_i \partial _{x^a}
        \end{split}
        \end{equation}
        \endproof

        \begin{thm}\label{Cp2}
        Let $\zt_i\colon E_i \rightarrow M,\ i=1,2,$ be vector
bundles over $M$ and let $\zC\colon E_1\rightarrow E_2$ be a
vector bundle morphism over the identity on $M$. Let $\zL_i$ be
a linear, 2-contravariant tensor on $E_i^\*$, $i=1,2$. Then
        $$ [\zC(X), \zC(Y)] ^{\zL_2} = \zC([X,Y]^{\zL_1})
                                                                        $$
        if and only if $\zL_2$ and $\zL_1$ are $\zC^\*$-related,
where $\zC^\* \colon E_2^\* \rightarrow  E_1^\*$ is the dual
morphism.
        \end{thm}
        \proof
         The equality $[\zC(X), \zC(Y)] ^{\zL_2} =
\zC([X,Y]^{\zL_1})$ is equivalent to the equality
        \begin{equation}\label{Fp6}
\{ (\zi_{E_1^\s*}X)\circ \zC^\*, (\zi_{E_1^\s*}Y) \circ \zC
^\*\}_{\zL_2} =
\{ \zi_{E_2^\s*}\zC(X), \zi_{E_2^\s*}\zC(Y)\}_{\zL_2} = \{
\zi_{E_1^\s*}X, \zi_{E_1^\s*}Y\}_{\zL_1}\circ \zC^\*.
        \end{equation}
        Since the exterior derivatives of functions $\zi_{E_1^\s*}X$
generate $\sT^\*E_1^\*$ over an open-dense subset ($E^\*_1$
minus the zero section), the equality~\ref{Fp6} holds if and only
if the tensors $\zL_1, \zL_2$ are $\zC^\*$-related.
        \endproof

        To the end of this section we assume that $\zL$ is
skew-symmetric, i.e., we consider pre-Lie algebroid structures only.

In this case, the bracket $[\,,\,]^\zL$, related to $\zL$, defined  on
sections of $\zt$ can be extended  in a standard way (cf.
\cite{G-U1,G-U2}) to the  graded bracket on $\zF(\zt)$.
We refer to this bracket as the Schouten-Nijenhuis bracket
and we denote it also by $[\,,\,]^\zL$.

        Moreover, we can define the 'exterior derivative' $\xd^\zL$
on $\zF(\zp)$ and the Lie derivative $\ll^\zL_X \colon \zF(\zp)
\rightarrow \zF(\zp)$ along a section $X\in \zG(M,E)$. Also the
Nijenhuis-Richardson bracket and the Fr\"olicher-Nijenhuis
bracket can be defined on $\zF_1(\zp) = \oplus_{n\in \Z}
\zF^n_1(\zp)$, where $\zF^n_1(\zp) = \zG(M,E\otimes \wedge ^n
E^\*)$. The definitions of these objects are analogous to the
definitions in the classical case (cf.~\cite{G-U2}).

        The bracket $[\,,\,]^\zL$ is a Lie bracket (or,
equivalently, $(\xd^\zL)^2 =0$) if and only if $\zL$ defines a
Lie algebroid structure, i.e., if and only if $\zL$ is a Poisson tensor. In
this case, all classical formulae   of differential
geometry, like $\ll^\zL_X \circ \xi_Y -\xi_Y \ll ^\zL_X =
\xi_{[X,Y]^\zL}$ etc., remain valid. We should also mention the
vertical tangent lift
        $$ \xv_\zt\colon \zG(M, \otimes_M ^kE) \rightarrow \zG(E,
\otimes _E^k\sT E)
                                                                $$
given, in local coordinates, by
        $$ \xv_\zt(f(x) e_{i_1}\otimes \cdots \otimes e_{i_k} ) = f(x)
\partial _{y^{i_1}} \otimes \cdots \otimes \partial _{y^{i_k}}.
                                                                $$
        In particular, $\xv_\zt(X\otimes Y) = \xv_\zt(X) \otimes
\xv_\zt(Y)$ (\cite{G-U2}). In the case of the tangent bundle, $E=\sT M$, the
vertical lift was denoted by $\vt$ in \cite{G-U1}. An analog
of the complete tangent lift $\dt$, studied for the tangent
bundle in  \cite{G-U1}, can be defined as follows.

        \begin{thm}\label{Cp3}
        Let $\zL$ be a linear bivector field on $E^\*$, which defines
a pre-Lie algebroid structure on a vector bundle $\zt\colon
E\rightarrow M$. Then, there exists a unique $\vt$-derivation of
order $0$
        $$ \dt^\zL \colon \otimes(\zt) \rightarrow \otimes(\ezt_E),
                                                                        $$
        which satisfies
        \begin{equation}\label{Fp7}
        \dt^\zL (f) = \zi_E \xd^\zL f \ \ \text{for}\ f\in C^\infty(M),
        \end{equation}
and
        \begin{equation}\label{Fp8}
        \zi_{\ssT^\s* E}(\dt^\zL X)\circ {\cal R} = \zi_{\ssT^\s*
E^\s*}([\zL, \zi_{E^\s*}X]) \ \ \text{for} \ X\in \zF^1(\zt),
        \end{equation}
        where $[\,,\,]$ is the Schouten bracket of multivector fields
on $E$ and ${\cal R}\colon \sT^\* E^\* \rightarrow \sT^\* E$ is the
canonical isomorphism of double vector bundles(see~\cite{M-X,U}).
Moreover, $\dt^\zL$ is a homomorphism of the Schouten-Nijenhuis
brackets:
        \begin{equation}\label{Fp9}
        \dt^\zL([X,Y]^\zL) = [\dt^\zL X, \dt^\zL Y],
        \end{equation}
if and only if $\zL$ is a Poisson tensor.
        \end{thm}
        {\sc Sketch of the proof.} Let us take  $X\in \zF^1(\zt)$.
The hamiltonian vector field $\G^\zL(X) = -[\zL,\zi_{E^\s*}X]$
is linear with respect to the tangent  vector bundle structure
$\sT\zt \colon \sT E^\* \rightarrow \sT M$ (\cite{U}). Hence,
the function $\zi_{\ssT^\s*E^\s*}[\zL, \zi_{E^\s*}X]$ is linear
with respect to both vector bundle structures on $\sT^\*E^\*$:
over $E$ and over $E^\*$.  It follows that there exists a unique
(linear) vector field $\dt^\zL X$ on $E$, such that
$\zi_{\ssT^\s*E^\s*}\G^\zL(X) = -(\zi_{\ssT^\s* E} \dt^\zL X) \circ
{\cal R}$.
        We have the formula
        $$ \dt^\zL(fX) = \dt^\zL(f)\xv_\zt(X) + \xv_\zt(f)\dt^\zL(X)
                                                                $$
and, consequently, we can extend $\dt^\zL$ to a
$\xv_\zt$-derivation on $\otimes(\zt)$.
        Finally, since $\cal R$ is an anti-Poisson isomorphism,
$\dt^\zL$ is a homomorphism of Schouten-Nijenhuis bracket if and
only if
        $$[\G^\zL(X), \G^\zL(Y)] = \G^\zL([X,Y]^\zL)$$
        for all $X, Y \in \zF^1(\zt)$, or, equivalently, if and only
if $\zL$ is a Poisson tensor.
        \hfill \vrule height4pt width6pt depth2pt

        {\bf Remark.} Let us define  a mapping
        $$ \I_E \colon  \zF^n_1(\zt) \rightarrow \zF^n(\ezt_E)
                                                        $$
by
        $$ \I_E(\zm\otimes X) = -\zi_E(\zm)\cdot \xv_\zt(X).
                                                $$
        It has been shown in \cite{G-U2} that $\I_E$ is a
homomorphism of the Nijenhuis-Richardson bracket into the
Schouten bracket. We have also a mapping
        $$ \G^\zL\colon \zF^n_1(\zp)\rightarrow \zF^n(\ezt_{E^\s*})
\colon K\mapsto \G^\zL(K) = [\zL,\I_{E^\s*}(K)],
                                                        $$
        which is, in the case of a Lie algebroid structure,  a
homomorphism of the Fr\"olicher-Nijenhuis bracket
$[\,,\,]^\zL_{F-N}$, associated to $\zL$, into the Schouten
bracket (\cite{G-U2}). The bracket $[\,,\,]^\zL_{F-N}$ is given
by the formula
                \begin{multline}
                [\zm\otimes X,\zn\otimes Y]^\zL = \zm\wedge \zn\otimes
[X,Y]^\zL + \zm\wedge \ll^\zL_X \zn \otimes Y -
\ll^\zL_Y \zm \wedge \zn \otimes X \\
                + (-1)^\zm (\xd^\zL\zm \wedge \xi_X\zn \otimes Y +
\xi_Y\zm \wedge \xd^\zL\zn \otimes X).
                        \end{multline}

\smallskip
        In local coordinates
        \begin{equation}\label{Fp10}
        \zL = \frac{1}{2}c^k_{ij}\zx_k \partial _{\zx_i}\wedge
\partial _{\zx_j} + \zd^a_i\partial _{\zx_i} \wedge \partial _{x^a}.
        \end{equation}
        Then,
        \begin{equation}\label{Fp11}
        \dt^\zL(f) = \frac{\partial f}{\partial x^a} \zd^a_j y^{j}
        \end{equation}
and
        \begin{equation}\label{Fp12}
        \dt^\zL(X^ie_{i}) = X^i\zd^a_i\partial _{x^a} + (X^i
c^k_{ji} + \frac{\partial X^k}{\partial x^a} \zd^a_j )y^j \partial _{y^k}.
        \end{equation}
It follows that, for $P = \frac{1}{2}P^{ij} e_i\wedge e_j $, we
have
        \begin{equation}\label{Fp13}
        \dt^\zL(P) = P^{ij}\zd^a_j \partial _{y^i}\wedge \partial
_{x^a} + (P^{kj} c^i_{lk} + \frac{1}{2}\frac{\partial
P^{ij}}{\partial x^a} \zd^a_l)y^l \partial _{y^i} \wedge
\partial _{y^j}.
        \end{equation}

        {\bf Remark.} For an arbitrary pseudo-Lie algebroid, we can
define the right and the left complete lifts with the use of the
right and the left hamiltonian vector fields instead of $[\zL,
\zi_E X]$.

        The following theorem describes the complete lifts in terms of
Lie derivatives and contractions.

        \begin{thm}\label{Cp4}
        Given a vector bundle $\zt\colon E\rightarrow M$ and a
linear bivector field $\zL$ on $E^\*$, we have
        \begin{description}
        \item[(a)] $\xv_\zt (X)(\zi_E\zm) = \xv_\zt(\xi_X\zm) =
\zt^\* \langle X, \zm\rangle $,
        \item[(b)]
        $ \dt^\zL (X)(\zi_E\zm) = \zi_E(\ll^\zL_X\zm)$,
        \end{description}
        where $X\in \zF^1(\zt)$ and $\zm\in \zF^1(\zp)$.
        \end{thm}
        \proof
        The part (a) has been proved in \cite{G-U2}, Theorem~15~c).
The part (b) follows from the following sequence of identities:
         $$\begin{aligned}
        \zp^\*_E (\dt^\zL(X) (\zi_E\zm))\circ {\cal R} &= \{
\zi_{\ssT^\s* E}(\dt^\zL X), \zp^\*_E(\zi_E\zm)\}_{\zL_E} \circ
{\cal R}\\
        &= \{ \zi_{\ssT^\s* E^\s*}(\xv_\zp \zm), \zi_{\ssT^\s*
E^\s*} ([\zL,\zi_{E^\s*}X])\}_{\zL_{E^*}} \\
        &= \zi_{\ssT^\s* E^\s*}[\xv_\zp\zm, [\zL,\zi_{E^\s*}X]] =
\zi_{\ssT^\s* E^\s*}[\G^\zL(X), \xv_\zp(\zm)] \\
        &= \zi_{\ssT^\s* E^\s*}(\xv_\zp(\ll^\zL_X\zm)) = \zp^\*_E(
\zi_E(\ll^\zL_X \zm)) \circ {\cal R},
                                        \end{aligned}$$
        where we used the equalities $[\G^\zL(X), \xv_\zp\zm] =
\xv_\zp (\ll^\zL_X\zm)$ (see \cite{G-U2}, Theorem~15~e)) and
$\zi_{\ssT^\s* E^\s*}\xv_\zp(\zm) = \zp^\*_E(\zi_E\zm) \circ
{\cal R}$.
        \endproof

        \begin{thm}\label{Cp5}
        If $P\in \zF^2(\zt)$, then $\dt^\zL(P)$ defines a pre-Lie
algebroid structure on $E^\*$ with the bracket
        \begin{equation}\label{Fp14}
        [\zm,\zn]^{\sdt^\zL(P)} = \ll_{P_\zm}^{\zL} \zn -
\ll_{P_\zn}^{\zL} \zm -\xd^{\zL}P(\zm,\zn),
        \end{equation}
        where $ P_\zm = \xi_\zm P$, and the anchor is given by
        $$ \za^{\sdt^\zL(P)}(\zm) = \za^\zL(P_\zm).
                                        $$
        \end{thm}
        \proof
        It is sufficient to consider $P= X\wedge Y$, $X,Y\in
\zF^1(\zt)$. Let us denote $\dt^\zL(X)$ and $\dt^\zL(Y)$ by
$\dot X $ and $\dot Y$, $\xv_\zt(X)$ and $\xv_\zt(Y)$ by $\bar
X$ and $\bar Y$,  $\ll^\zL$ and $\xd^\zL$ by $\ll$ and $\xd$\,.
Then, we have
         $$\begin{aligned}
                \{ \zi_E\zm, \zi_E \zn\}_{\sdt^\zL(P)} &= \dot
X(\zi_E\zm) \bar Y (\zi_E\zn ) - \dot X(\zi_E\zn) \bar Y
(\zi_E\zm ) - \bar X(\zi_E\zm) \dot Y (\zi_E\zn ) + \bar
X(\zi_E\zn) \dot Y (\zi_E\zm )\\
                &= \zi_E\left(\langle Y,\zn\rangle \ll_X(\zm) - \langle
Y,\zm\rangle \ll_X(\zn) - \langle X,\zn\rangle \ll_Y(\zm) +
\langle X,\zm\rangle \ll_X(\zn)\right)\\
                &= \zi_E\left(\ll_{P_\zm}\zn - \ll_{P_\zn}\zm -\xd(\langle
X,\zm\rangle \langle Y,\zn\rangle - \langle X,\zn\rangle \langle
Y,\zm\rangle)\right),
                \end{aligned}$$
where we used Theorem~\ref{Cp4}. Now, we have
        $$
        [\zm,f\zn]^{\sdt^\zL(P)}= \ll_{P_\zm}(f\zn) - \ll_{P_{f\zn}}\zm
-\xd P(\zm,f\zn) = f[\zm,\zn]^{\sdt^\zL(P)} + (\ll_{P_\zm}f)\zn,
                        $$
        so that $\za^{\sdt^\zL(P)}(\zm) = \ll_{P_\zm}(f) =
\za^\zL(P_\zm)(f)$.
        \endproof

        In the case of the canonical Lie algebroid  on the
tangent bundle $\ezt_M \colon \sT M\rightarrow M$, associated
with the canonical Poisson tensor $\zL_M$ on $\sT^\* M$, our
definition of $\dt^\zL(P)$ gives the standard tangent complete
lift $\dt$. Moreover, the bracket \ref{Fp14} of 1-forms is the
bracket introduced independently in \cite{Fu,Ko,KS-M} and
corresponding to the lift $\dt P$ (cf. \cite{Co,G-U1}).

        {\bf Example.}
        Let us consider a Lie algebroid over a point, i.e., a  real Lie
algebra ${\frak g}$ with a basis $e_1,\dots,e_m$, and its dual
space  ${\frak g}^\*$ with the dual basis $e^{*1},\dots,e^{*m}$.
We have also the corresponding (linear) coordinate system
$\zx_1,\dots, \zx_m $ on ${\frak g}^\*$ and the coordinate system
$y^1,\dots, y^m$  on ${\frak g}$. The linear Poisson structure
$\zL$ on ${\frak g}^\*$, associated with the Lie bracket
$[\,,\,]^\zL$ on ${\frak g}$, is the well-known
Kostant-Kirillov-Souriau tensor
        $$ \zL = \frac{1}{2} c^k_{ij}\zx_k \partial
_{\zx_i}\wedge \partial _{\zx_j}.
                                $$
        Here $c^k_{ij}$ are the structure constants with respect to
the chosen basis. The exterior derivative $\xd^\zL \colon
\wedge {\frak g}^\* \rightarrow  \wedge {\frak g}^\* $ is the dual
mapping to the Lie bracket:
        $$ \xd^\zL(\zm)(X,Y) = \langle \zm, [Y,X]^\zL\rangle ,
                                $$
i.e., $\xd^\zL$ is the Chevalley cohomology operator. For $x \in
\zF^1(\zt) = {\frak g} $, the tangent complete lift $\dt^\zL(x)$
is the fundamental vector field of the adjoint representation,
corresponding to $x$:
                $$ \dt^\zL(e_i) = c^k_{ji}y^j \partial _{y^k}.
                                $$

\section{\hspace{-4mm}.\hspace{2mm}NIJENHUIS TENSORS AND POISSON-NIJENHUIS
STRUCTURES FOR LIE ALGEBROIDS}
        Let a vector bundle $\zt \colon E \rightarrow M$ be given a
pseudo-Lie algebroid structure, associated with a tensor field
$\zL$ on $E^\*$, and let $\widetilde{N}\colon E\rightarrow E $
be a vector bundle morphism over the identity. We can represent
$\widetilde{N}$, as well as its dual ${\widetilde{N}}^\*$, by a
tensor field $N\in \zF^1_1(\zp)$. This tensor field defines
operations in $\zF^1(\zt)$ and $\zF^1(\zp)$, which we denote by
the same symbol $\xi_N$. If $N = X_i\otimes \zm^i,\ X_i\in
\zF^1(\zt),\ \zm^i \in \zF^1(\zp) $, the operation $\xi_N$ is
given by the formulae
        $$   \xi_N X = \langle X,\zm^i\rangle X_i \ \ \text{and}\ \
\xi_N \zm = \langle X_i,\zm \rangle \zm^i,
                                $$
        where $X\in \zF^1(\zt),\ \zm \in \zF^1(\zp)$.
        In the notation of \cite{KS-M}, $\xi_N X = NX$ and $\xi_N
\zm =\!{\phantom{|}}^t N\zm$. It is obvious that we can extend
$\xi_N$ to a derivation of the tensor algebra, putting
        \begin{equation}\label{Fp15}
        \xi_N(A\otimes B) = (\xi_N A)\otimes B + A\otimes (\xi_N B).
        \end{equation}
        Using $N$, we can deform the bracket $[\,,\,]^\zL$ to a
bracket $[\,,\,]^\zL_N$ on $\zF^1(\zt)$ by the formula
        \begin{equation}\label{Fp16}
        [X,Y]^\zL_N = [NX,Y]^\zL +[X,NY]^\zL -N[X,Y]^\zL.
        \end{equation}

        \begin{thm}\label{Cp6}
        The deformed bracket \ref{Fp16} defines on $E$ a pseudo-Lie
algebroid structure, with the anchors $(\za^\zL_N)_l = \za^\zL_l
\circ \widetilde{N}$ and $(\za^\zL_N)_r = \za^\zL_r \circ \widetilde{N}$. The
associated tensor field is given by
        $$
        \zL_N = \ll_{\I_{E^*}(N)}\zL,
                                        $$
        where $ \ll_{\I_{E^*}(N)}$ is the standard Lie derivative
along the vector field $\I_{E^\s*}(N)$.

        If $\zL$ is skew-symmetric, then $\zL_N$ is also
skew-symmetric and the Schouten-Nijenhuis bracket, induced by
$\zL_N$, can be written in the form, similar to \ref{Fp16},
        \begin{equation}\label{Fp17}
        [X,Y]^{\zL_N}= [X,Y]^\zL_N = [\xi_N X,Y]^\zL + [X,\xi_N
Y]^\zL -\xi_N([X,Y]^\zL)
        \end{equation}
        for $X,Y \in \zF(\zt)$. Moreover,
        \begin{equation}\label{Fp18}
        \xd^{\zL_N} = \xi_N \circ \xd^\zL - \xd^\zL \circ \xi_N.
        \end{equation}
        \end{thm}
        The proof is based on the following Lemma.
        \begin{lem}
        For $X \in\zF^1(\zt)$, we have
        \begin{equation}\label{Fp19}
        \zi_{E^\s*}(NX) = - \ll_{\I_{E^*}(N)}(\zi_{E^\s*})
        \end{equation}
        and
        \begin{equation} \label{Fp20}
        \ll_{\I_{E^*}(N)} \xv_\zp (\zm) = \xv_\zp(\xi_N \zm)
        \end{equation}
for $X\in \zF^1(\zt),\ \zm \in\zF^1(\zp)$.
        \end{lem}
        \proof
        Let $N = X_i\otimes \zm^i$, $X_i \in \zF^1(\zt)$ and $\zm^i \in
\zF^1(\zp)$. We have
        $$
        \begin{aligned}
        \zi_{E^\s*}(NX) &= \zi_{E^\s*}(\langle X,\zm^i\rangle X_i) =
\zp^\* (\langle X,\zm^i\rangle ) \zi_{E^\s*}(X_i)\\
        &= \xv_\zp (\langle X,\zm^i\rangle ) \zi_{E^\s*}(X_i).
        \end{aligned}
                                $$
        On the other hand,
        $$
        \begin{aligned}
        -\ll_{\I_{E^*}(N)}(\zi_{E^\s*}X) &= (\zi_{E^\s*}(X_i)\xv_\zp(
\zm^i)) (\zi_{E^\s*}X) = \zi_{E^\s*}(X_i)\xv_\zp(\zm^i)(\zi_{E^\s*}X) \\
        &= \zi_{E^\s*}(X_i) \xv_\zp(\langle X, \zm^i\rangle ),
        \end{aligned}
                                $$
        according to Theorem~15~c) in \cite{G-U2}.

        Similarly,
        $$
        \begin{aligned}
        [\I_{E^*}(N),\xv_\zp(\zm)] &= [\zi_{E^\s*}(X_i) \xv_\zp(\zm_i),
\xv_\zp(\zm)] \\
        &= - \xv_\zp (\zm^i) \wedge [\zi_{E^*}(X_i),\xv_\zp(\zm)],
        \end{aligned}
                                $$
        since the vertical vector fields commute. Following
Theorem~15~c) in \cite{G-U2}, we get
        $$ [\zi_{E^\s*}(X_i),\xv_\zp(\zm)] = -\xv_\zp (\xi_{X_i}\zm)
                        $$
                and, consequently,
        $$[\I_{E^*}(N), \xv_\zp(\zm)] = \xv_\zp(\zm^i \wedge
\xi_{X_i}\zm) = \xv_\zp(\xi_N \zm).
                        $$
        \endproof

        {\sc Proof of Theorem~\ref{Cp6}.}
        Using Lemma and properties of the Lie derivative, we get
        $$
        \begin{aligned}
        \zi_{E^\s*}([X,Y]^\zL_N) &= \zi_{E^\s*}\left( [NX,Y]^\zL
 + [X,NY]^\zL -N[X,Y]^\zL \right)\\
        &= - \{ \ll_{\I_{E^*}(N)}(\zi_{E^\s*}X), \zi_{E^\s*}Y\}_\zL
- \{\zi_{E^\s*}X, \ll_{\I_{E^*}(N)} (\zi_{E^\s*}Y)\}_\zL +\\
        & \qquad +\ll_{\I_{E^*}(N)}\{\zi_{E^\s*}X, \zi_{E^\s*}Y \}_\zL\\
        &= \{ \zi_{E^\s*}X, \zi_{E^\s*}Y\}_{\ll_{\I_{E^*}(N)}\zL}.
                \end{aligned}
                                $$
        The general form \ref{Fp17} of the corresponding Schouten
bracket follows inductively from the Leibniz rule for the Schouten
bracket $[\,,\,]^\zL$ and from \Ref{Fp15}.

        In order to prove \Ref{Fp18}, we, again, use Lemma and
\cite{G-U2} Theorem~15~d):
        $$
        \begin{aligned}
        \xv_\zp(\xd^{\zL_N}\zm) &= [\zL_N, \xv_\zp\zm] = [[\I_{E^\s*}N,
\zL],\xv_\zp \zm]\\
        &= [\I_{E^\s*}N, [\zL,\xv_\zp\zm ]] - [\zL,[\I_{E^\s*}N,
\xv_\zp\zm]] \\
        &= [\I_{E^\s*}N, \xv_\zp(\xd_\zL \zm)] -[\zL, \xv_\zp
(\xi_N \zm)] = \xv_\zp(\xi_N \xd^\zL \zm - \xd^\zL \xi_N \zm).
        \end{aligned}
                                $$
        \hfill \vrule height4pt width6pt depth2pt

        In local coordinates, we have
        $$
        \begin{aligned}
          N &= N^i_j e_i\otimes e^{*j},\\
        \zL &= c^k_{ij}\zx_k \partial _{\zx_i} \otimes \partial
_{\zx_j} + \zd^a_i \partial _{\zx_i} \otimes \partial _{x^a} -
\zs^a_i \partial _{x^a}\otimes \partial _{\zx_i},\\
        \I_{E^\s*}N &= N^i_k \zx_i \partial _{\zx_k},
        \end{aligned}
                                $$
        and
        \begin{multline}
        \zL_N = \left( c^k_{lj} N^l_i + c^k_{il} N^l_j - c^l_{ij}
N^k_l + \zd^a_i \frac{\partial N^k_j}{\partial x^a} -  \zs^a_j
\frac{\partial N^k_i}{\partial x^a} \right) \zx_k \partial _{\zx_i}
\otimes \partial _{\zx_j}\\
        + N^l_i\zd^a_l \partial _{\zx_i}\otimes \partial _{x^a}
        - N^l_i\zs^a_l \partial _{x^a}\otimes \partial _{\zx_i}.
        \end{multline}

        \begin{thm}\label{Cp7}
        For $X \in \otimes (\zt)$ and skew-symmetric $\zL$, we have
        $$ \dt^{\zL_N} (X) = \dt^\zL(\xi_N X) -
\ll_{\I_{E}(N)}\dt^\zL (X).
                        $$
        \end{thm}
        \proof
        Since $\dt^{\zL_N}$ and $\dt^\zL$ are $\xv_\zt$-derivations
of order 0 on $\otimes (\zt)$ and $\ll_{\I_{E}(N)} \xv_\zt(X) =0$
($\I_{E}(N)$ is vertical), it is enough to consider the case
$X\in \zF^1(\zt)$. For such $X$
        $$
        \begin{aligned}
        \zi_{\ssT^\s* E} \left(\dt^{\zL_N} (X)\right)\circ {\cal R}
&=  \zi_{\ssT^\s* E^\s*} [\zL_N, \zi_{E^\s*}X] = \zi_{\ssT^\s*
E^\s*} [\ll_{\I_{E^*}(N)}\zL, \zi_{E^\s*}X]\\
        &= \zi_{\ssT^\s* E^\s*}\left(\ll_{\I_{E^*}(N)} [\zL,
\zi_{E^\s*}X \right)  - \zi_{\ssT^\s* E^\s*}[\zL,
\ll_{\I_{E^*}(N)}(\zi_{E^\s*}X)] .
        \end{aligned}
                                $$
         Since $\ll_{\I_{E^*}(N)} = -\zi_{E^\s*}(NX)$ (\Ref{Fp19}), then
        $$ -\zi_{\ssT^\s* E^\s*}[\zL,
\ll_{\I_{E^*}(N)}(\zi_{E^\s*}X)] = \zi_{\ssT^\s* E^\s*}[\zL,
\zi_{E^\s*}(NX)] = \zi_{E^\s*}\left( \dt^\zL(NX)\right) \circ
{\cal R}.
                        $$
        On the other hand,
        $$
        \begin{aligned}
        \zi_{\ssT^\s* E^\s*}\left(\ll_{\I_{E^*}(N)} [\zL, \zi_{
E^\s*}X \right) &= \{\zi_{\ssT^\s* E^\s*}(\I_{E^*}(N)),
\zi_{\ssT^\s* E^\s*}[\zL, \zi_{E^\s*}X] \}_{\zL_{E^\s*}}\\
        &= \{\zi_{\ssT^\s*E}(\I_EN)\circ {\cal R},
\zi_{\ssT^\s*E}(\dt^\zL(X))\circ {\cal R} \}_{\zL_{E^\s*}}\\
        &= -\{\zi_{\ssT^\s*E}(\I_EN), \zi_{\ssT^\s* E}(\dt^\zL (X)) \}_{\zL_E}
\circ {\cal R}\\
        &= -\left(\zi_{\ssT^\s* E}[\I_EN, \dt^\zL(X)]\right)
\circ {\cal R}\\
        &= -\zi_{\ssT^\s* E} \left(\ll_{\I_{E}(N)}\dt^\zL(X)\right),
        \end{aligned}
                                $$
where we used the equality $\zi_{\ssT^\s* E}(\I_{E}N)
=\zi_{\ssT^\s* E^\s*}(\I_{E^*}(N))$. Since ${\cal R}$ is an
isomorphism and $\zi_{\ssT^\s* E}$ is injective, the theorem follows.
        \endproof

\medskip
        In local coordinates, for $\zL$ as in \Ref{Fp10}, we have
        \begin{multline}
         \dt^{\zL_N} (X^ie_i) = X^iN^k_i \zd^a_k \partial _{x^a} +\\
        +\left( X^i ( N^k_j c^n_{ki} + N^k_i c^n_{jk} -N^n_k c^k_{ji}
+ \zd^a_j \frac{\partial N^n_i}{\partial x_a} - \zd^a_i
\frac{\partial N^n_j}{\partial x_a}) + \frac{\partial
X^n}{\partial {x^a}} N^k_j \zd^a_k \right) y_j \partial _{y^n}.
                          \end{multline}

        {\bf Remark.}  If we treat the Schouten brackets $B^\zL =
[\,,\,]^\zL$ and $B^\zL_N = [\,,\,]^\zL_N$ as bilinear operators
on $\zF(\zt)$, then  fomula~\Ref{Fp16} means
        \begin{equation} B^\zL_N = [\xi_N, B^\zL]_{N-R},
        \end{equation}
where $[\,,\,]_{N-R}$ is the Nijenhuis-Richardson bracket of
multilinear graded operators of a graded space in the sense of
\cite{LMS}. Similarly, \Ref{Fp18} means that
        \begin{equation}\label{Fp21}
        \dt^{\zL_N} = [\xi_N, \dt^\zL]_{N-R}.
        \end{equation}

This interpretation will be used later, together with
the Jacobi identity for the $[\,,\,]_{N-R}$.

        \begin{defn}\label{Cp8}
        A tensor $N\in \zG(M,E\otimes E^\*)$ is called a {\em
Nijenhuis tensor} for $\zL$ (or, for a Lie algebroid structure
defined by $\zL$), if the {\em Nijenhuis torsion}
        \begin{equation} T^\zL_N(X,Y) = N[X,Y]^\zL_N - [NX,NY]^\zL
        \end{equation}
        vanishes for all $X,Y\in \zG(E).$
        \end{defn}

        The classical version of the following  is well known
(cf. \cite{KS-M}).

        \begin{thm}\label{Cp9} \makebox[1cm]{}

        \begin{description}
        \item[(a)] $N$ is a Nijenhuis tensor for $\zL$ if and only if
$\zL$ and $\zL_N = \ll_{\I_{E^*}(N)}\zL$ are $\widetilde{N}^\*$-related.
        \item[(b)] The Nijenhuis torsion corresponds to the
Fr\"olicher-Nijenhuis bracket:
        $$ T^\zL_N(X,Y)  = \frac{1}{2}[N,N]^\zL_{F-N}
(X,Y).
                                        $$
        \item[(c)]
        $$  [[B^\zL,\xi_N]_{R-N}, \xi_N]_{N-R} = 2T^\zL_N + [B^\zL,
\xi_{N^2}]_{N-R},
                        $$
        where $(X_i\otimes \zm^i)^2 = \langle X_i,\zm^j \rangle
\zm^i \otimes X_j$.
        \item[(d)]   If $N$ is a Nijenhuis tensor, then $\zL_N$ is a
Poisson tensor.
        \end{description}
        \end{thm}
        \proof

        (a) \ \ Since $\zL_N = \ll_{\I_{E^*}(N)}\zL$ induces the
deformed bracket $B^\zL_N = [\,,\,]^\zL_N$, this part follows
from Theorem~\Ref{Cp2}.

        (b) \ \  Let $N=X_i\otimes \zm^i$, then
        $$ [NX,NY]^\zL = \langle X,\zm^i\rangle \langle Y,
\zm^j\rangle [X_i, X_j]^\zL +   \langle X,\zm^i\rangle
\ll^\zL_{X_i} (\langle Y,\zm^j\rangle) X_j -  \langle Y,\zm^j\rangle
\ll^\zL_{X_j} (\langle X,\zm^i\rangle) X_i
                        $$
and
        \begin{multline}
        N[X,Y]^\zL_N = N\left( \langle X,\zm^i\rangle [X_i,Y]^\zL -
\ll^\zL_{Y}(\langle X,\zm^i\rangle)X_i + \langle Y,\zm^j\rangle
[X,X_j] \right.\\
        + \left.\ll^\zL_{X}(\langle Y,\zm^j\rangle) X_j - \langle
[X,Y]^\zL ,\zm^i\rangle X_i\right)\\
        = \langle X,\zm^i\rangle\langle [X_i,Y]^\zL,\zm^j\rangle X_j
- \ll^\zL_{Y} (\langle X,\zm^i\rangle) \langle X_i,\zm^j\rangle
X_j\\
        + \langle Y,\zm^j\rangle\langle [X,X_j]^\zL,\zm^i\rangle X_i
+ \ll^\zL_{X} (\langle Y,\zm^j\rangle) \langle X_j,\zm^i\rangle
X_i - \langle [X,Y]^\zL,\zm^i\rangle \langle X_i,\zm^j\rangle X_j.
        \end{multline}

        Hence, using properties of Lie derivatives, we get
        $$\begin{aligned}
        T^\zL_N (X,Y) &= \langle X,\zm^i\rangle \langle
Y,\zm^j\rangle [X_i,X_j]^\zL + \langle X,\zm^i\rangle \langle Y,
\ll^\zL_{X_i}\zm^j \rangle X_j \\
        &\  - \langle Y,\zm^j\rangle \langle X, \ll^\zL_{X_j}\zm^i
\rangle X_i + \xd \zm^i (X,Y) \langle X_i,\zm^j\rangle X_j \\
        &= \left( \frac{1}{2} \zm^i\wedge \zm^j \otimes [X_i,
X_j]^\zL + \zm_i\wedge \ll^\zL_{X_i}\zm^j\otimes X_j  +\xd^\zL\zm^i
\wedge \xi_{X_i} \zm^j \otimes X_j \right) (X,Y)\\
        &= \frac{1}{2} [N,N]^\zL_{F-N} (X,Y).
        \end{aligned}$$

        (c) \ \
        \begin{multline}
        [[B^\zL, \xi_N]_{R-N}, \xi_N]_{R-N} (X,Y)=\\
        =([N^2X,Y]^\zL + [NX,NY]^\zL  -N [NX,Y]^\zL) + ([NX, NY]^\zL
+ [X, N^2Y]^\zL - N[X, NY]^\zL) \\
        - (N [NX,Y]^\zL + N [X, NY]^\zL - N^2[X,Y]^\zL)\\
        = 2([NX,NY]^\zL - N([NX,Y]^\zL + [X,NY]^\zL - N[X,Y]^\zL)) \\
        + [N^2X,Y]^\zL +[X, N^2Y]^\zL - N^2[X,Y]^\zL\\
        =2 T^\zL_N(X,Y) +[B^\zL,\xi_{N^2}]_{N-R}(X,Y).
        \end{multline}

        (d) \ \ The Schouten bracket induced by $\zL_N$ is given by
$[\xi_N,B^\zL ]_{N-R}$ and it is known from general theory
\cite{LMS}, that it defines a graded Lie algebra structure if and
only if its Nijenhuis-Richardson square vanishes. Using the
graded Jacobi identity for $[\,,\,]_{N-R}$, we get
        $$
        \begin{aligned}
        [[\xi_N,B^\zL]_{N-R}, [\xi_N,B^\zL]_{N-R}]_{N-R} &=
[[[\xi_N,B^\zL]_{N-R} ,\xi_N]_{N-R},  B^\zL]_{N-R}\\
        &= -2[T^\zL_N, B^\zL]_{N-R} +[[\xi_{N^2},B^\zL]_{N-R},
B^\zL]_{N-R}\\
        &=0,
        \end{aligned}
                        $$
        since $T^\zL_N =0$ and  $[B^\zL,B^\zL]_{N-R} =0$
($[\,,\,]^\zL$ is a Lie bracket) implies that
$\left(\ad^{N-R}_{B^\zL}\right)^2 =0$.

        \endproof

        The following theorem is, essentially, due to Mackenzie and
Xu (\cite{M-X}).

        \begin{thm}\label{Cp10}
                Let $\zL$ be a Poisson tensor on $E^\*$ and let $P\in
\zF^2(\zt)$. Then
        \begin{description}
        \item[(a)] $\dt^\zL(P)$ induces a pre-Lie algebroid structure
on $E^\*$, with the bracket and the anchor described in
Theorem~\Ref{Cp5}. The exterior derivative, induced by
$\dt^\zL(P)$, is given by
        \begin{equation}\label{Fp22}
        \xd^{\sdt^\zL(P)}(X) = [P,X]^\zL.
        \end{equation}
        Moreover,
        \begin{equation}\label{Fp23}
        \frac{1}{2}[P,P]^\zL(\zm,\zn,\zg) =\langle
\widetilde{P}({[\zm,\zn]}^{\dt^\zL(P)}) -[P_\zm,P_\zn]^\zL, \zg \rangle
        \end{equation}
for all $\zm,\zn,\zg \in \zF^1(\zp)$ and $P$ is a Poisson tensor
for $\zL$ (i.e., $[P,P]^\zL =0$) if and only if $\zL$ and
$\dt^\zL(P)$ are $-\widetilde{P}$-related, where
$\widetilde{P}(\zm) =P_\zm = \xi_\zm P$.
        \item[(b)] if $P$ is, in addition, a Poisson tensor for
$\zL$, then $\dt^\zL(P)$ is a Poisson tensor and Poisson tensors
$\zL, \dt^\zL(P)$ induce a Lie bialgebroid structure on bundles
$E$ and $E^\*$, i.e.,
        \begin{equation}\label{Fp24}
        \xd^\zL\left( [\zm,\zn]^{\sdt^\zL(P)}\right) = [\xd^\zL \zm,
\zn]^{\sdt^\zL(P)} +(-1)^{\zm+1}[\zm, \xd^\zL \zn]^{\sdt^\zL(P)}.
        \end{equation}
        \end{description}
        \end{thm}

        \proof
        The proof of \Ref{Fp23} is completely analogous to the proof
in the classical case (see \cite{KS-M}). The remaining part of
(a) follows from Theorem~\Ref{Cp2}. Part(b) is proved in
\cite{M-X}.

        \endproof

        {\bf Remark.} Due to the result of Kosmann-Schwarzbach
(\cite{KS}), \Ref{Fp24} is equivalent to
        \begin{equation}\label{Fp25}
        [P,[X,Y]^\zL]^\zL = [[P,X]^\zL,Y]^\zL + (-1)^X[X,[P,Y]^\zL]^\zL,
        \end{equation}
        which is a special case of the graded Jacobi identity for the
bracket $[\,,\,]^\zL$.

        The fact that $\dt^\zL(P)$ is a Poisson tensor,
if $[P,P]^\zL=0$, is a direct consequence of \Ref{Fp10}. The
converse to this is not true, in general, as  shows the
following example.

        {\bf Example 2.} For a Lie algebroid over a point, i.e., for
a Lie algebra ${\frak g}$ with the bracket $[\,,\,]^\zL$,
corresponding to a Kirillov-Kostant-Souriau tensor $\zL$ on
${\frak g}^\*$,  $P\in \wedge ^2{\frak g}$ is a Poisson tensor
for $\zL$ if and only if $P$  is  an r-matrix satisfying the
classical Yang-Baxter equation $[P,P]^\zL = 0$.

        On the other hand, $\dt^\zL(P)$ is a Poisson tensor if and
only if $\dt^\zL([P,P]^\zL) =0$ which means, that $\ad_\zx
[P,P]^\zL =0$ for all $\zx\in {\frak g}$, i.e., the equation
 $\dt^\zL([P,P]^\zL) =0$ is the modified Yang-Baxter equation.

        \begin{defn}\label{Cp11}
Let $P\in \zF^2(\zt)$ be a Poisson tensor with respect to a Lie
algebroid structure on $E$, associated to a Poisson tensor $\zL$
on $E^\*$, and let $N \in \zF^1_1(\zt)$ be a Nijenhuis tensor for
$\zL$. We call the pair $(P,N)$ a {\em Poisson-Nijenhuis
structure} for $\zL$ if the following two conditions are
satisfied:
        \begin{enumerate}
        \item $NP =PN^\*$, where $NP(\zm,\zn) = P(\zm, \xi_N\zn)$  and
$PN^\*(\zm,\zn) = P(\xi_N \zm, \zn)$,
        \item $\dt^{\zL_N}(P) = \left(\dt^\zL(P)\right)_N. $
        \end{enumerate}
        \end{defn}

        {\bf Remark.} Since $NP +PN^\* =\xi_N P$ and, according to
Theorems \Ref{Cp6} and \Ref{Cp7},
                $$ \left(\dt^\zL(P)\right)_N = \ll_{\I_{E}(N)} \dt^\zL(P),
                        $$
        $$ \dt^{\zL_N} (P) = \dt^\zL(\xi_N P) - \ll_{\I_{E}(N)} \dt^\zL(P),
                        $$
the condition (2) can be replaced by

                (2') $ \ll_{\I_{E}(N)} \dt^\zL(P) =
\left(\dt^\zL(P)\right)_N =\dt^\zL(NP) $.

\medskip
                \begin{thm}\label{Cp12}
        If $(P,N)$ is a Poisson-Nijenhuis structure for $\zL$ then
$NP$ is a Poisson tensor for $\zL$ and we have the following
commutative diagram of Poisson mappings between Poisson
manifolds.
        $$\begin{CD}
        (E^\*,\zL)  @> -\widetilde{P}>> (E,\dt^\zL(P)) \\
        @VV \widetilde{N}^\*V   @VV\widetilde{N}V \\
        (E^\*, \zL_N) @>-\widetilde{P}>> (E,\dt^\zL(NP) =
\left(\dt^\zL(P)\right)_N)
                        \end{CD},$$
where $\zL_N=\ll_{\I_{E^*}(N)} \zL$ and
$\left(\dt^\zL(P)\right)_N = \ll_{\I_{E}(N)}\dt^\zL(P)$.
Moreover, every structure from the left-hand side of this diagram
constitutes a Lie bialgebroid structure with every right-hand side
structure.
        \end{thm}
        \proof
        The tensors $\zL_N$ and $\dt^\zL(P)$ are Poisson. The
mappings $-\widetilde{P}\colon (E^\*,\zL) \rightarrow (E,
\dt^\zL(P))$ and $\widetilde{N}^\* \colon (E,\zL) \rightarrow
(E^\*, \zL_N)$ are Poisson, in view of Theorems \Ref{Cp9} and
\Ref{Cp10}. The assumption $NP = PN^\*$ implies that the diagram
is commutative. To show that the mapping $-\widetilde{P} \colon
(E^\*, \zL_N) \rightarrow (E,\left(\dt^\zL(P)\right)_N)$ is
Poisson, it is enough to check that, under the assumption  $NP =
PN^\*$, the vector fields $\I_{E^*}(N)$ and $\I_{E}(N)$ are
$-\widetilde{P}$-related. One can do it easily.
        Since $\zL$ and $\dt^\zL(P)$ are $-\widetilde{P}$-related,
also $\zL_N = [\I_{E^*}(N),\zL]$ and $\left(\dt^\zL(P)\right)_N=
[\I_{E}(N), \dt^\zL(P)]$ are $-\widetilde{P}$-related. Hence, the
equality $\left(\dt^\zL(P)\right)_N = \dt^\zL(NP)$ implies that
$\zL$ and $\dt^\zL(NP)$ are $-\widetilde{NP}$-related and,
according to Theorem \Ref{Cp10}~a), $NP$ is a Poisson tensor for
$\zL$.

The fact that the mapping $\widetilde{N}\colon (E, \dt^\zL(P)) \rightarrow
(E,\left(\dt^\zL(P)\right)_N)$ is Poisson follows from the
identity
        \begin{multline}
        \langle X, [N,N]^{\sdt^\zL(P)}_{F-N}(\za,\zb)\rangle  = \langle
[N,N]^\zL_{F-N} (X, P_\zb),\za\rangle \\
        + 2\langle X, C^\zL(P,N)(\xi_N\za,\zb)\rangle  - 2\langle
NX, C^\zL(P,N)(\za,\zb)\rangle,
        \end{multline}
        where $C^\zL(P,N)(\za,\zb) = [\za,\zb]^{\sdt^\zL(NP)} -
[\za,\zb]^{\sdt^\zL(P)}_N$. This is a generalization of an
analogous identity in \cite{KS-M}, with a completely parallel
proof.

        The pairs $(\zL,\dt^\zL(P))$ and $(\zL, \dt^\zL(NP))$
constitute Lie bialgebroids by Theorem~\Ref{Cp10}~b), since $P$
and $NP$ are Poisson tensors for $\zL$.

        Similarly, $(\zL_N, \dt^\zL(NP)) = \dt^{\zL_N}(P)$ constitute
a Lie bialgebroid, since $P$ is a Poisson tensor for $\zL_N$
($\zL_N$ and $\dt^{\zL_N}(P)$ are $-\widetilde{P}$-related).

        To show that the pair $(\zL_N, \dt^\zL(P))$ forms a Lie
bialgebroid, we have to prove that $\xd^{\zL_N}=\xd^\zL_N$ is a
derivation of the Schouten bracket $B =[\,,\,]^{\sdt^\zL(P)}$,
i.e., $[\xd^\zL_N,B]_{N-R} =0$. Since, due to~\Ref{Fp18}, $\xd_N^\zL
= [\xi_N, \xd^\zL]_{N-R}$ and $[\xd^\zL,B]_{N-R} =0$, we get
        \begin{multline}
        [\xd_N^\zL,B]_{N-R} = [[\xi_N,\xd^\zL]_{N-R},B]_{N-R} =
[\xi_N,[\xd^\zL ,B]_{N-R}]_{N-R}\\
        + [\xd^\zL, [\xi_N,B]_{N-R}]_{N-R} = [\xd^\zL,
B_{(\sdt^\zL(P))_N}]_{N-R} = 0,
        \end{multline}
in view of the fact that $[\xi_N,B]_{N-R}$ is the bracket
associated to $\left(\dt^\zL(P)\right)_N$, for which $\xd^\zL$ is a
derivation.
    \phantom{aa}    \endproof

        {\bf Remark.} The above diagram is a dualization of a similar
diagram in \cite{KS-M}.

        In the case of the canonical Lie algebroid on $E=\sT M$, the
fact that $((\zL_M)_N, \dt P)$ constitutes a Lie bialgebroid is
equivalent to the fact that $(P,N)$ is a Poisson-Nijenhuis
structure, as it was recently shown in \cite{KS}. This is due to
the formulae
        \begin{align}
        A(f,g) &= \langle (NP-PN^\*)\xd^\zL g, \xd^\zL f \rangle ,\\
        A(\xd^\zL f,g)& = C^\zL(P,N)(\xd^\zL f, \xd^\zL g) + \xd^\zL A(f,g),
        \end{align}
        where $A = [\xd^{\zL_N} ,B_{\sdt^\zL(P)}]_{N-R}$,
and the fact that $A$ satisfies a Leibniz rule and $\xd^\zL f$
generate $\sT^\* M$, for $\zL =\zL_M$.

        In general, this is not true and we can have $((\zL)_N,
\dt P)$ being a Lie bialgebroid with $(P,N)$ not being
Poisson-Nijenhuis structure for $\zL$, even if we assume the equality  $NP
=PN^\*$, as shows the following example.

        {\bf Example.} As a Lie algebroid over a single point, let us
take a Lie algebra ${\frak g}$ spanned by $\zx_1,\zx_2,\zx_3, \zx_4$
with the bracket defined by $\zL =\zx_3 \partial _{\zx_1} \wedge
\partial _{\zx_2}$. The tensor  $P=\partial _{\zx_2}\wedge
\partial _{\zx_4}$ is a Poisson tensor with $\dt^\zL P = y_1
\partial _{y_3} \wedge \partial _{y_4}$.

        The tensor
        $$N= -\zx_1\otimes y_1 + \sum_{i=2}^4 \zx_i\otimes y_i$$
is a Nijenhuis tensor for $\zL$ and $\zL_N =-\zL.$
        Moreover, $NP =PN^\* =P$, so that $(\zL, \dt^\zL(NP))$
constitutes a Lie bialgebroid. In this case, however,
$\dt^{\zL_N}P = -\dt^\zL P = -\dt^\zL(NP)$ and (PN) is not a
Poisson-Nijenhuis structure.

        It is easy to see that, as in the classical case, a
Poisson-Nijenhuis structure for a Lie algebroid induces a whole
hierarchy of compatible Poisson structures and Nijenhuis tensors
(see \cite{KS-M}). Since this theory goes quite parallel to the
classical case, we will not present details here.

  \end{document}